\documentclass{amsart}

\usepackage{amsfonts}
\usepackage{amsmath}
\usepackage{amsthm}
\usepackage{amssymb}
\usepackage[english]{babel}
\usepackage{graphics}
\usepackage{latexsym}
\usepackage{longtable} \usepackage{mathrsfs}
\usepackage{graphicx}
\usepackage{wrapfig} 
\usepackage[pdftex,colorlinks=true,linkcolor=blue,citecolor=cyan]{hyperref}

\newtheorem{theorem}{Theorem}
\newtheorem{corollary}[theorem]{Corollary}
\newtheorem{lemma}[theorem]{Lemma}

\newtheorem{claim}[theorem]{Claim}
\newtheorem{example}[theorem]{Example}

\theoremstyle{definition}
\newtheorem{definition}[theorem]{Definition}

\newcommand{\mH}{\mathcal{H}}

\newcommand{\R}{\mathbb{R}}
\newcommand{\N}{\mathbb{N}}

\newcommand{\D}{\mathrm{D}}

\newcommand{\noi}{\noindent}
\newcommand{\ms}{\medskip}

\newcommand{\al}{\alpha}

\newcommand{\De}{\Delta}
\newcommand{\e}{\varepsilon}

\newcommand{\Om}{\Omega}
\newcommand{\om}{\omega}

\newcommand{\larrow}{\longrightarrow}

\newcommand{\lmapsto}{\longmapsto}
\newcommand{\ri}{\rightarrow}

\newcommand{\p}{\partial}
\newcommand{\sub}{\subseteq}

\newcommand{\by}{\times}

\newcommand{\Lip}{\mathrm{Lip}}

\newcommand{\ess}{\mathrm{ess}}
\newcommand{\diam}{\mathrm{diam}}
\newcommand{\dist}{\mathrm{dist}}

\newcommand{\bt}{\begin{theorem}}\newcommand{\et}{\end{theorem}}
\newcommand{\bd}{\begin{definition}}\newcommand{\ed}{\end{definition}}
\newcommand{\bl}{\begin{lemma}}\newcommand{\el}{\end{lemma}}
\newcommand{\beq}{\begin{equation}}\newcommand{\eeq}{\end{equation}}
\newcommand{\bc}{\begin{claim}}\newcommand{\ec}{\end{claim}}
\newcommand{\bex}{\begin{example}}\newcommand{\eex}{\end{example}}
\newcommand{\bcor}{\begin{corollary}}\newcommand{\ecor}{\end{corollary}}
\newcommand{\bp}{\begin{proof}}\newcommand{\ep}{\end{proof}}

\numberwithin{equation}{section}

\begin{document}

\title[Global $W^{1,p}$ bounds for Harmonic functions]{A Remark on Global $W^{1,p}$ Bounds for Harmonic Functions with Lipschitz Boundary Values}

\author{Nikos Katzourakis}
\address{Department of Mathematics and Statistics, University of Reading, Whiteknights, PO Box 220, Reading RG6 6AX, Reading, UK}
\email{n.katzourakis@reading.ac.uk}


\date{}


\begin{abstract} In this note we show that gradient of Harmonic functions on a smooth domain with Lipschitz boundary values is pointwise bounded by a universal function which is in $L^p$ for all finite $p\geq 1$. 
\end{abstract}

\maketitle


\section{Introduction} \label{section1}

Kellogg in \cite{K} pioneered the study of the boundary behaviour of the gradient of Harmonic functions on a bounded domain. Roughly speaking, he established that in a domain of $\R^3$ near a boundary region which can represented as the graph of a planar function, the gradient of any Harmonic function is continuous up to the boundary provided that the gradient of the boundary function and of the Harmonic function are Dini continuous themselves on the boundary. The celebrated theory of Schauder estimates \cite{GT} establishes strong relevant results for general uniformly elliptic PDEs, providing interior and global H\"older bounds for solutions and their derivatives in terms of the H\"older norms of the boundary values of the solution and the right hand side of the PDE. The Schauder theory has been improved and extended by many authors, but typically for second order elliptic PDEs with boundary values of the solutions and right hand sides of the PDEs in the H\"older spaces $C^{2,\al}$ or $C^{1,\al}$, in order to obtain uniform estimates for the solutions in the respective H\"older spaces. 

In \cite{GH} Gilbarg-H\"ormander have extended Schauder theory to include hypotheses of lower regularity of the boundary values of the solution, of the boundary of the domain and of the coefficients of the equations. Troianiello \cite{T} relaxed further some conditions of Gilbarg-H\"ormander \cite{GH}. In the paper \cite{HS} Hile-Stanoyevitch, extending an older result of Hardy-Littlewood \cite{HL}, proved that the gradient of a Harmonic function with Lipschitz continuous boundary values is pointwise bounded up to a constant by the logarithm of a multiple of the inverse of the distance to the boundary.

However, it appears that in none of these results, even for the special case of the Laplacian, there is an explicit \emph{global bound in $L^p$} for the gradient of Harmonic functions which have just Lipschitz boundary values and not $C^{1,\al}$. In this note establish the following consequence of the result of Hile-Stanoyevitch:

\begin{theorem} \label{theorem1} Let $n\geq 2$, $\Om\sub \R^n$ a bounded open set with $C^2$ boundary. Let also $g :\p\Om \ri \R$ with $g \in \Lip(\p\Om)$, that is $g\in C^0(\p\Om)$ and 
\[
\Lip(g,\p\Om)\,:=\, \sup_{x,y\in \p\Om,\, x\neq y} \frac{|g(x)-g(y)|}{|x-y|}\, <\,\infty.
\]

\begin{enumerate}
\item \label{(i)} There exists a positive function $f_{\Om,n} : \Om \ri (0,\infty)$ depending on $\Om,n$ such that 
\beq \label{1.1}
\ \ f_{\Om,n} \, \in \bigcap_{p\in[1,\infty)}\! L^p\Om) \cap C^0(\Om)
\eeq
and if $u \in C^2(\Om)\cap C^0(\overline{\Om})$ is the Harmonic function solving
\beq  \label{1.2}
\ \  \left\{
\begin{array}{rl}
\De u\, =\, 0, & \text{in }\Om, \smallskip
\\
u\,=\,g, & \text{on }\p\Om,
\end{array}
\right.
\eeq
then we have the estimate
\beq  \label{1.3}
\ \ \big| \D u(x)\big|\, \leq\, \Lip(g,\p\Om)\, f_{\Om,n}(x), \quad x\in\Om.
\eeq

\ms

\item \label{(ii)} Let $(g^m)_1^\infty \sub \Lip(\p\Om)$ satisfy for some $C>0$
\beq  \label{1.4}
\ \  \Lip(g^m,\p\Om)\, +\, \max_{\p\Om} |g^m| \, \leq\,C, \quad m\in \N.
\eeq
Let also $(u^m)_1^\infty \in C^2(\Om)\cap C^0(\overline{\Om})$ be the Harmonic functions solving 
\beq  \label{1.5}
\ \ \left\{
\begin{array}{rl}
\De u^m\, =\, 0,\ \ & \text{in }\Om, \smallskip
\\
u^m\,=\,g^m, & \text{on }\p\Om.
\end{array}
\right.
\eeq
Then, $(u^m)_1^\infty$ is strongly precompact in $\bigcap_{p=1}^\infty W^{1,p}(\Om)$ and if 
\beq  \label{1.6}
\ \ \ \text{$g^{m_k} \larrow g$ \ in $C^0(\overline{\Om})$, as $k\ri \infty$,} 
\eeq
then there is a unique limit point $u \in C^2(\Om)\cap C^0(\overline{\Om})$ of the subsequence $(u^{m_k})_1^\infty$ such that along perhaps a further subsequence
\beq  \label{1.7}
\ \ \ \ u^{m_k} \larrow u \ \ \text{ in }W^{1,p}(\Om) \ \forall \, p\geq 1, \text{ as }k\ri\infty,
\eeq
and the limit function $u$ solves
\[
\ \ \left\{
\begin{array}{rl}
\De u\, =\, 0, & \text{in }\Om, \smallskip
\\
u\,=\,g, & \text{on }\p\Om.
\end{array}
\right.
\]
\end{enumerate}

\end{theorem}

The motivation to derive the above integrability result and its consequences comes from certain recent advances in generalised solutions of nonlinear PDE and vectorial Calculus of Variations in the space $L^\infty$ (\cite{Ka4} and \cite{Ka2, Ka3}). The vectorial counterparts of Harmonic functions provide useful energy comparison maps since they are ``stable" in $L^p$ for all $1<p<\infty$.
\ms

\section{Proofs}

Our notation is either self-explanatory or otherwise standard as e.g.\ in \cite{E}, \cite{Ka}. The starting point of our proof is the following estimate of Hile-Stanoyevitch: under the hypotheses of Theorem \ref{theorem1}, the gradient $\D u$ of a Harmonic function $u\in C^2(\Om)\cap C^0(\overline{\Om})$ which solves \eqref{1.2} with $g\in \Lip(\p\Om)$ satisfies the logarithmic estimate
\beq  \label{2.1}
\big| \D u(x)\big| \,\leq\, C(\Om,n)\, \Lip(g,\p\Om)\, \ln \left( \frac{\diam(\Om)}{\dist(x,\p\Om)}\right) ,\ \ \ x\in \Om.
\eeq
for some $C$ depending just on $\Om$ (and the dimension). In \eqref{2.1}, $\diam(\Om)$ is the diameter of the domain and $\dist(x,\p\Om)$ the distance of $x$ from the boundary:
\[
\begin{split}
\diam(\Om)\, &:=\, \sup\big\{|x-y|\, :\, x,y\in \Om \big\},
\\
\dist(x,\p\Om)\, &:=\, \inf\big\{|x-z|\, :\, z\in \p \Om \big\}.
\end{split}
\]

\ms

\noi \textbf{Proof of \eqref{(i)} of Theorem \ref{theorem1}.} Fix $\e >0$ smaller than the diameter of $\Om$ and consider the inner open $\e$ neighbourhood of $\Om$:
\[
\Om^\e\,:=\, \big\{x\in \Om\,:\, \dist(x,\p\Om)>\e \big\}.
\]
It is well known that (see e.g.\ \cite{GT})
\[
\dist(\cdot,\p \Om)\, \in \, W^{1,\infty}_{\text{loc}}(\R^n)
\]
and 
\beq \label{2.2}
\big|\D\, \dist(\cdot,\p \Om)\big|\, =\, 1, \ \ \ \text{a.e.\ on }\Om.
\eeq
Let $p\in[1,\infty)$. By the Co-Area formula (see e.g.\ [\cite{EG}, Proposition 3, p.\ 118]) applied to the function
\[
\R^n \ni \, x \lmapsto \chi_{\Om^\e}(x) \left(\ln \left( \frac{\diam(\Om)}{\dist(x,\p\Om)}\right)\right)^p \, \in \R
\]
(where $\chi_{\Om^\e}$ is the characteristic function of ${\Om^\e}$), we have 
\beq \label{2.3}
\begin{split}
 \int_{\Om^\e} & \left(\ln \left( \frac{\diam(\Om)}{\dist(x,\p\Om)}\right)\right)^p dx
\, = 
\\
&= \, \int_\e^{\diam(\Om)} \left(\int_{\{\dist(\cdot,\p\Om)=t\}}  \frac{\left(\ln \left( \dfrac{\diam(\Om)}{\dist(z,\p\Om)}\right)\right)^p}{\big|\D\, \dist(z,\p \Om)\big|}   d\mH^{n-1}(z)\right) dt
\end{split}
\eeq
where $\mH^{n-1}$ is the $(n-1)$-dimensional Hausdorff measure. By using \eqref{2.2}, \eqref{2.3} simplifies to
\[
\begin{split}
 \int_{\Om^\e} & \left(\ln \left( \frac{\diam(\Om)}{\dist(x,\p\Om)}\right)\right)^p dx
\, = 
\\
&= \, \int_\e^{\diam(\Om)} \left(\int_{\{\dist(\cdot,\p\Om)=t\}}  \left(\ln \left( \dfrac{\diam(\Om)}{\dist(z,\p\Om)}\right)\right)^p d\mH^{n-1}(z)\right) dt
\end{split}
\]
Further, since 
\[
\text{$\dist(z,\p\Om) = t$,\ for all $z\in \{\dist(\cdot,\p\Om)=t\}$,}
\] 
by setting
\beq
\label{2.4}
I^{\e,p}\,:=\,  \int_{\Om^\e} \left(\ln \left( \frac{\diam(\Om)}{\dist(x,\p\Om)}\right)\right)^p dx
\eeq
we obtain
\beq \label{2.4a}
\begin{split}
I^{\e,p} \, & = \, \int_\e^{\diam(\Om)} \left(\int_{\{\dist(\cdot,\p\Om)=t\}}  \left(\ln \left( \frac{\diam(\Om)}{t}\right)\right)^p   d\mH^{n-1}(z)\right) dt
\\
& = \,  \int_\e^{\diam(\Om)} \left(\ln \left( \frac{\diam(\Om)}{t}\right)\right)^p \mH^{n-1}\Big(\{\dist(\cdot,\p\Om)=t\}\Big)\,  dt.
\end{split}
\eeq
As a consequence of the regularity of the boundary, standard results regarding the equivalence between the Hausdorff measure and the Minkowski content for rectifiable sets (see e.g.\ [\cite{AFP}, Section 2.13, Theorem 2.106]) imply that there is a $C=C(\Om)$ such that
\[
\underset{0<t<\diam(\Om)}{\ess\,\sup} \, \mH^{n-1}\Big(\{\dist(\cdot,\p\Om)=t\}\Big) \,\leq\, C(\Om)
\]
and hence the inequality \eqref{2.4a} gives
\beq \label{2.5}
I^{\e,p} \, \leq\,C(\Om) \int_\e^{\diam(\Om)} \left(\ln \left( \frac{\diam(\Om)}{t}\right)\right)^p   dt.
\eeq
By the change of variables
\[
\om\,:=\, \frac{\diam(\Om)}{t}
\]
we can rewrite the estimate \eqref{2.5} as
\[
I^{\e,p} \, \leq\,C(\Om)\, \diam(\Om) \int_1^{\diam(\Om)/\e}   \frac{(\ln \om)^p}{\om^2}\,   d\om
\]
and by enlarging perhaps the constant $C(\Om)$, we rewrite this as
\beq \label{2.6}
I^{\e,p} \, \leq\,C(\Om)  \int_1^{\diam(\Om)/\e}  \left( \frac{\ln \om }{\om^{2/p}}\right)^p   d\om .
\eeq

\begin{claim} \label{claim2} We have that
\[
\lim_{\e\ri 0}\, I^{\e,p}\, \leq \, C(\Om,n,p)\,<\infty.
\]
\end{claim}

\noi {\bf First proof of Claim } \ref{claim2} (proposed by one of the referees): By using the following known property of the Gamma function
\[
\int_1^\infty \frac{\ln^p x}{x^2}\, dx\, =\, \Gamma(1+p)
\]
we readily conclude. \qed

\ms

\noi {\bf Second proof of Claim } \ref{claim2}: We now give a direct argument without quoting special functions. Consider now the function
\[
g(\om)\,:=\, \frac{\ln \om }{\om^{2/p}}, \ \ \ \ g\, :\  (1,\infty)\ri (0,\infty).
\]
Since
\[
g'(\om)\,=\, \frac{\, 1 - (2/p)\ln \om \,}{\om^{(2/p) +1}}, 
\]
we have that $g$ is strictly increasing on $(1,e^{2/p})$ and strictly decreasing on $(e^{2/p},\infty)$. Further, note that $t \mapsto g^p(t)$ also enjoys the exact same monotonicity properties since $s\mapsto s^p$ is strictly increasing. 
Moreover, since 
\[
e^{2/p}\, \leq\, 10
\]
for all $p\in [1,\infty)$ and by using that $\e \mapsto I^{\e,p}$ is decreasing (in view of \eqref{2.4}), we have
\[
\begin{split}
\lim_{\e\ri 0} I^{\e,p}\, & \leq\, C(\Om) \left[ \int_1^{10}\frac{(\ln \om)^p}{\om^2}\, dt\ +\ \int_{10}^\infty \frac{(\ln \om)^p}{\om^2}\, dt \right]
\\
 & \leq\, C(\Om) \left[ 10 \left(\sup_{1<\om<10}\frac{(\ln \om)^p}{\om^2}\right)\ +\ \int_{10}^\infty \frac{(\ln \om)^p}{\om^2}\, dt \right]
 \\
 & =\, C(\Om) \left[ 10 \left( \frac{(\ln \om)^p}{\om^2}\right) \! \Big|_{\om=e^{2/p}}\ +\ \int_{10}^\infty \frac{(\ln \om)^p}{\om^2}\, dt \right]
\end{split}
\]
and hence
\[
\begin{split}
\lim_{\e\ri 0} I^{\e,p}
\,  & \leq \, C(\Om) \left[ 10  \frac{(2/p)^p}{e^{4/p}} \ +\ \sum_{k=10}^\infty\int_k^{k+1}\frac{(\ln \om)^p}{\om^2}\, dt \right]
\\
 & \leq\, C(\Om) \left[ 10  \frac{(2/p)^p}{e^{4/p}} \ +\ \sum_{k=10}^\infty \left(\sup_{k<\om<k+1}\frac{(\ln \om)^p}{\om^2}\right) \right]
\end{split}
\]
which gives
\beq \label{2.7}
\lim_{\e\ri 0} I^{\e,p}\, \leq\, C(\Om) \left[ 10  \frac{(2/p)^p}{e^{4/p}} \ +\ \sum_{k=10}^\infty \frac{(\ln k)^p}{k^2}\right].
\eeq
Now we show that the series
\[
S\,:=\, \sum_{k=10}^\infty \frac{(\ln k)^p}{k^2}
\]
converges.

\ms
\noi \textbf{Method 1}: Since the sequence 
\[
\frac{(\ln k)^p}{k^2}, \ \ \ \ k\,=\,10,11,12,...
\]
is decreasing, by the Cauchy condensation test the series $S$ converges if and only if 
\[
\sum_{m=10}^\infty  2^m\, A_{2^m} <\, \infty, \ \ \ \ \ A_k\,:=\, \frac{(\ln k)^p}{k^2}.
\]
Since
\[
\frac{2^{m+1}\, A_{2^{m+1}}}{2^m\, A_{2^m}}\, =\, \frac{(m+1)^p\, (\ln 2)^p \,2^{-m-1}}{m^p\, (\ln 2)^p\, 2^{-m}}\,=\, \frac{1}{2} \left(1+ \frac{1}{m}\right)^p \larrow \,\frac{1}{2},
\]
as $m\ri \infty$, by the Ratio test we have that 
\[
\sum_{k=10}^\infty \frac{(\ln k)^p}{k^2}\, =\, C(p)\, <\,\infty 
\]
since $\sum_{m=10}^\infty  2^m\, A_{2^m}$ converges. 

\ms

\ms
\noi \textbf{Method 2} (proposed by one of the referees): By repeated applications of the del Hospital rule ($p\in \N$), we have
\[
\lim_{k\ri \infty} \, \dfrac{\, \dfrac{(\ln k)^p}{k^2} \,}{ \dfrac{1}{k^{3/2}} }\, = \, \lim_{k\ri \infty} \, \frac{p!}{k}\, =\, 0
\]
and hence 
\[
 \dfrac{(\ln k)^p}{k^2}\, \leq\, \dfrac{1}{k^{3/2}}
\]
for $k\in \N$ large enough. Since the series
\[
\sum_{k=10}^\infty\, \dfrac{1}{k^{3/2}}
\]
converges and hence so does $S$ by the comparison test. 
\ms

In either case, by \eqref{2.4} and \eqref{2.7} we have that there is a constant $C(\Om,n,p)$ depending only on $\Om,n,p$ such that
\beq \label{2.8}
\begin{split}
\int_{\Om} \left(\ln \left( \frac{\diam(\Om)}{\dist(x,\p\Om)}\right)\right)^p dx \, &=\, \lim_{\e\ri 0}I^{\e,p}
\\ 
& \leq\, C(\Om,n,p).
\end{split}
\eeq
By combining \eqref{2.8} with \eqref{2.1}, we see that by setting
\[
f_{\Om,n}(x) \, :=\, C(\Om,n)\ln \left( \frac{\diam(\Om)}{\dist(x,\p\Om)}\right),\ \ \ x\in\Om
\]
\eqref{(i)} of Theorem \ref{theorem1} is established.   \qed 
\ms

\ms

\noi \textbf{Proof of \eqref{(ii)} of Theorem \ref{theorem1}.} Let $u^m$ solve \eqref{1.5}. By standard interior bounds on the derivatives of Harmonic functions in terms of their boundary values (see e.g.\ \cite{GT}) and \eqref{1.4}, we have that the Hessians $(\D^2u^m)_1^\infty$ are bounded in $C^0(\Om,\R^{n\by n})$, that is uniformly over the compact subsets of $\Om$. The same is true for the $3$rd order derivatives as well; thus, for any $\Om'\Subset \Om$, there is $C(\Om')$ such that 
\[
\sum_{k=1}^3 \, \big\|\D^ku^m\big\|_{C^0(\Om')}\, \leq\, C(\Om') \, \|u^m\|_{C^0(\Om)}
\]
and by the Maximum Principle we have
\[
\|u^m\|_{C^0(\Om)}\, \leq\, \max_{\p\Om}|g^m|\, \leq\,C.
\]
As a consequence, 
\[
\left.
\begin{array}{c}
\big| \D^ku^m(x)-\D^ku^m(y)\big| \,\leq\, C(\Om')|x-y|, \ms\\ 
 \ \ \ \ \ \ \, \big|\D^ku^m(x)\big|\,\leq\, C(\Om'), 
\end{array}
\right\}\ \ \ \  x,y\in \Om',\ k=0, 1,2,3,\ m\in \N
\]
and by the Ascoli-Arzela theorem, the sequence 
\[
\Big(u^m,\D u^m,D^2u^m\Big)_{m=1}^\infty 
\]
is precompact uniformly over the compact subsets of $\Om$. Again by \eqref{1.4}, we have 
\[
\left.
\begin{array}{c}
\big| g^m(x)-g^m(y)\big| \,\leq\, C|x-y|,\  \ \ \ms\\ 
\big|g^m(x)\big|\,\leq\, C, 
\end{array}
\right\}\ \ \ \  x,y\in \p\Om,\ m\in \N
\]
which gives that $(g^m)_1^\infty$ is bounded and equicontinuous on $\p\Om$. Thus, by the Ascoli-Arzela theorem and by the lower semicontinuity of the Lipschitz seminorm with respect to uniform convergence, there is a subsequence $(g^{m_k})_1^\infty$ and $g\in \Lip(\p\Om)$ such that
\[
g^{m_k} \larrow g,\ \ \ \text{as }k\ri \infty \text{ in }C^0(\p \Om).
\]
Along perhaps a further subsequence, by the above bounds on $(u^m)_1^\infty \sub C^2(\Om)\cap C^0 (\overline{\Om})$, there is $u\in C^2(\Om)$ such that
\beq \label{2.9}
\left\{\ \
\begin{split}
u^{m_k} & \larrow u, \ \ \ \ \ \text{ in }C^0(\Om), 
\\
\D u^{m_k} & \larrow \D u, \ \ \ \text{ in }C^0(\Om,\R^n), 
\\
\D^2 u^{m_k}& \larrow \D^2u, \ \ \text{ in }C^0(\Om,\R^{n\by n}), 
\end{split}
\right.
\eeq
as $k\ri \infty$. By passing to the limit in the equation $\De u^m=0$ we get that $\De u=0$. Since the measure of $\Om$ is finite, for any $p\in [1,\infty)$ by H\"older inequality we have that
\[
\|u^m\|_{L^p(\Om)}\, \leq\, |\Om|^{1/p}\Big( \|u^m\|_{C^0(\Om)} \Big) \leq\, C(\Om,p).
\]
By item \eqref{(i)} of the theorem and by \eqref{1.4}, we have that 
\[
\|\D u^m\|_{L^p(\Om)}\, \leq\, C(\Om,n,p).
\]
Hence, we have the bound
\[
\| u^m\|_{W^{1,p}(\Om)}\, \leq\, C(\Om,n,p), \ \ \ p\geq 1.
\]
By the Morrey embedding theorem, by choosing $p>n$ we have that $(u^{m_k})_1^\infty$ is precompact in $C^0(\overline{\Om})$ and hence by \eqref{2.9} we have that 
\beq \label{2.10}
u^{m_k} \larrow u, \ \ \ \text{ in }C^0(\overline{\Om})\ \text{ as }k\ri\infty.
\eeq
Hence, $u=g$ on $\p\Om$ and as a consequence $u$ solves the limit Dirichlet problem. Finally, if $E\sub \Om$ is a measurable subset, by the H\"older inequality we have that
\[
\begin{split}
\int_E \big|Du^m(x)\big|^p\,dx\, & \leq\, |E|^{1-\frac{p}{p+1}}\left(\int_E \big|Du^m(x)\big|^{p+1}\,dx\right)^{\!\frac{p}{p+1}}
\\
& =\, |E|^{1-\frac{p}{p+1}}\Big(\|\D u^m\|_{L^{p+1}(\Om)}\Big)^p
\\
& \leq\, |E|^{1-\frac{p}{p+1}}C(\Om,n,p).
\end{split}
\]
Hence, the sequence of gradients $(\D u^{m_k})_1^\infty$ is $p$-equi-integrable on $\Om$. By \eqref{2.9}, we have 
\[
\text{$\D u^{m_k} \larrow \D u$\ \ in measure on $\Om$, \ as $k\ri \infty$.}
\] 
Since $\Om$ has finite measure, the Vitali Convergence theorem (e.g.\ \cite{FL}) implies that 
\[
\text{$\D u^{m_k} \larrow \D u$\ \ in $L^p(\Om)$, \ as $k\ri \infty$.}
\] 
Item \eqref{(ii)} of Theorem \ref{theorem1} has been established.    \qed 
\ms

\ms

\noi {\bf Acknowledgement.} The author would like to thank the referees of this paper most warmly for the careful reading of the manuscript and for providing thought-through alternative simpler proofs of certain of the original arguments.
\ms

\ms

\end{document}